\documentclass[reqno]{amsart}
\usepackage{amsmath}
\usepackage{amsthm}
\usepackage{amsfonts}
\usepackage{amssymb}
\usepackage{latexsym}
\usepackage{amscd}
\usepackage{stmaryrd}
\usepackage{color}

\definecolor{darkgreen}{rgb}{0.0625,0.64,0.0625}
\usepackage[
            pdfstartview=FitH,
            CJKbookmarks=true,
            bookmarksnumbered=true,
            bookmarksopen=true,
            colorlinks,
            pdfborder=001,
            linkcolor=blue,
            anchorcolor=green,
            citecolor=red
            ]{hyperref}

\theoremstyle{plain}
\newtheorem{thm}{Theorem}[section]
\newtheorem{cor}[thm]{Corollary}
\newtheorem{lem}[thm]{Lemma}

\numberwithin{equation}{section}

\newfont{\scyr}{wncyr10 scaled 550}
\def\shuffle{\,\mbox{\bf \scyr X}\,}

\def\proof{\noindent {\bf Proof.\;}}

\allowdisplaybreaks
\begin{document}

\title{Stuffle product formulas of multiple zeta values}

\date{\today\thanks{The first author is supported by the National Natural Science Foundation of
China (Grant No. 11471245) and Shanghai Natural Science Foundation (grant no. 14ZR1443500). The authors would like to thank the referee for helpful suggestions, which improve the paper greatly.} }

\author{Zhonghua Li \quad and \quad Chen Qin}

\address{School of Mathematical Sciences, Tongji University, No. 1239 Siping Road,
Shanghai 200092, China}

\email{zhonghua\_li@tongji.edu.cn}

\address{School of Mathematical Sciences, Tongji University, No. 1239 Siping Road,
Shanghai 200092, China}

\email{2014chen\_qin@tongji.edu.cn}



\begin{abstract}
Using the combinatorial descriptions of stuffle product, we obtain recursive formulas for the stuffle product of multiple zeta values and of multiple zeta-star values. Then we apply the formulas to prove several stuffle product formulas with one or two strings of $z_p$'s. We also describe how to use our formulas in general cases.
\end{abstract}

\maketitle


\section{Introduction and statement of main results}\label{Sec:Intro}

For integers $k_1\geqslant2$ and $k_2,\ldots,k_n\geqslant1$, the multiple zeta value and the multiple zeta-star value are defined by
\begin{align*}
&\zeta(k_1,k_2,\ldots,k_n)=\sum\limits_{m_1>m_2>\cdots>m_n>0}\frac{1}{m_1^{k_1}m_2^{k_2}\cdots m_n^{k_n}},\\
&\zeta^\star(k_1,k_2,\ldots,k_n)=\sum\limits_{m_1\geqslant m_2\geqslant \cdots\geqslant m_n\geqslant 1}\frac{1}{m_1^{k_1}m_2^{k_2}\cdots m_n^{k_n}},
\end{align*}
respectively.
In case of $n=1$, we get Riemann zeta values in both cases, which are special values of the Riemann zeta function at positive integer arguments.

Using the definitions we can express a product of two multiple zeta values (multiple zeta-star values) as a sum of multiple zeta values (multiple zeta-star values). For example, for multiple zeta values we have
\begin{align*}
\zeta(2)\zeta(3)&=\sum\limits_{m>0}\frac{1}{m^2}\times\sum\limits_{n>0}\frac{1}{n^3}=\left(\sum\limits_{m>n>0}
+\sum\limits_{n>m>0}+\sum\limits_{m=n>0}\right)\frac{1}{m^2n^3}\\
&=\zeta(2,3)+\zeta(3,2)+\zeta(5),
\end{align*}
and for multiple zeta-star values, we have
\begin{align*}
\zeta^\star(2)\zeta^\star(3)&=\sum\limits_{m\geqslant 1}\frac{1}{m^2}\times\sum\limits_{n\geqslant 1}\frac{1}{n^3}=\left(\sum\limits_{m\geqslant n\geqslant 1}+\sum\limits_{n\geqslant m\geqslant 1}-\sum\limits_{m=n>0}\right)\frac{1}{m^2n^3}\\
&=\zeta^\star(2,3)+\zeta^{\star}(3,2)-\zeta^\star(5).
\end{align*}
Such products are called stuffle products (or harmonic shuffle products).

The algebraic settings are given in \cite{Hoffman} for multiple zeta values, and in \cite{Shuichi Muneta} for multiple zeta-star values. Since we only treat the stuffle product here, we give the following algebraic descriptions. Let $Y=\{z_k\mid k=1,2,\ldots\}$ be an alphabet with noncommutative letters, and let $Y^{\ast}$ be the set of all words generated by letters in $Y$, which contains the empty word $1$. Let $\mathfrak{h}^1=\mathbb{Q}\langle Y\rangle$ be the noncommutative polynomial algebra over the rational number field $\mathbb{Q}$ generated by the alphabet $Y$, which is spanned by $Y^{\ast}$ as a rational vector space. And let $\mathfrak{h}^0$ be the subalgebra of $\mathfrak{h}^1$ defined by
$$\mathfrak{h}^0=\mathbb{Q}+\sum\limits_{k=2}^\infty z_{k}\mathfrak{h}^1.$$
Corresponding to stuffle products of multiple zeta values and of multiple zeta-star values, we define two bilinear commutative products $\ast$ and $\bar{\ast}$ on $\mathfrak{h}^1$ by the rules
\begin{align*}
&1\ast w=w\ast 1=w,\\
&z_kw_1\ast z_lw_2=z_k(w_1\ast z_lw_2)+z_l(z_kw_1\ast w_2)+z_{k+l}(w_1\ast w_2);\\
&1 \,\bar{\ast}\, w=w\,\bar{\ast}\, 1=w,\\
&z_kw_1\,\bar{\ast}\, z_lw_2=z_k(w_1\,\bar{\ast}\, z_lw_2)+z_l(z_kw_1\,\bar{\ast}\, w_2)-z_{k+l}(w_1\,\bar{\ast}\, w_2),
\end{align*}
where $w,w_1,w_2\in Y^{\ast}$ and $k,l$ are positive integers. Under each of these two products, $\mathfrak{h}^1$ becomes a commutative algebra, and $\mathfrak{h}^0$ is still a subalgebra of $\mathfrak{h}^1$.

To relate the algebra $\mathfrak{h}^0$ to multiple zeta values and to multiple zeta-star values, we define  rational linear maps $Z: \mathfrak{h}^0 \rightarrow \mathbb{R}$ and $Z^{\star}: \mathfrak{h}^0\rightarrow \mathbb{R}$ by $Z(1)=Z^{\star}(1)=1$ and
$$Z(z_{k_1}z_{k_2}\cdots z_{k_n})=\zeta(k_1,k_2,\ldots,k_n),\quad Z^{\star}(z_{k_1}z_{k_2}\cdots z_{k_n})=\zeta^{\star}(k_1,k_2,\ldots,k_n),$$
where $n,k_1,k_2,\ldots,k_n$ are positive integers with $k_1>1$, and $\mathbb{R}$ is the field of real numbers.
Then it was shown in \cite{Hoffman} and in \cite{Shuichi Muneta} that $Z:(\mathfrak{h}^0,\ast)\rightarrow\mathbb{R}$ and $Z^{\star}:(\mathfrak{h}^0,\bar{\ast})\rightarrow\mathbb{R}$ are algebra homomorphisms, respectively. More precisely, for any $w_1,w_2\in \mathfrak{h}^0$, we have
\begin{align*}
&Z(w_1\ast w_2)=Z(w_1)Z(w_2),\\
&Z^{\star}(w_1\,\bar{\ast}\, w_2)=Z^\star(w_1)Z^\star(w_2).
\end{align*}
Hence to get stuffle product formulas, we can work on $\mathfrak{h}^0$ first. For example, since
\begin{align*}
&z_2\ast z_3=z_2z_3+z_3z_2+z_5,\\
&z_2\,\bar{\ast}\, z_3=z_2z_3+z_3z_2-z_5,
\end{align*}
applying the map $Z$ and $Z^{\star}$ respectively, we get
$$\zeta(2)\zeta(3)=\zeta(2,3)+\zeta(3,2)+\zeta(5), \quad \zeta^\star(2)\zeta^\star(3)=\zeta^\star(2,3)+\zeta^{\star}(3,2)-\zeta(5)$$
as proved above.

There is another commutative product $\shuffle$ among multiple zeta values, which is deduced from the iterated integral representation of multiple zeta values and is called the shuffle product (see \cite{Hoffman,Zagier} for details). During the study of multiple zeta values, several shuffle product formulas were found. For example, as generalizations of Euler's decomposition formula, Eie and Wei essentially studied the shuffle product $z_kz_1^m\shuffle z_lz_1^n$ in \cite{Eie-Wei}, Lei, Guo and Ma obtained the formula for shuffle product $z_kz_1^m\shuffle z_{l_1}z_1^{n_1}z_{l_2}z_1^{n_2}$ in \cite{Lei-Guo-Ma}, and we studied the shuffle products of words with more strings of $z_1$'s in \cite{Li-Qin}. For the stuffle product, Chen obtained a formula for the product $z_kz_p^m\ast z_lz_p^n$ in \cite{Chen}. Chen expressed this product by simple products of the form $z_p^i\ast z_p^j$. We find that one can treat general products with more strings of $z_p$'s such as
\begin{align}
z_{k_1}z_p^{m_1}z_{k_2}z_p^{m_2}\cdots z_{k_r}z_p^{m_r}\ast z_{l_1}z_p^{n_1}z_{l_2}z_p^{n_2}\cdots z_{l_s}z_p^{n_s}
\label{Eq:Stu-MZV}
\end{align}
and
\begin{align}
z_{k_1}z_p^{m_1}z_{k_2}z_p^{m_2}\cdots z_{k_r}z_p^{m_r}\,\bar{\ast}\, z_{l_1}z_p^{n_1}z_{l_2}z_p^{n_2}\cdots z_{l_s}z_p^{n_s}.
\label{Eq:Stu-MZSV}
\end{align}
In fact, the general products can be expressed by simple products of the forms $z_p^i\ast z_p^j$ and $z_p^i\,\bar{\ast}\,z_p^j$, respectively. For the simple products, we have the following result.

\begin{thm}\label{Thm:Simple-Product}
For a positive integer $p$ and two nonnegative integers $a,b$, we denote by $G_b^a$ the set of all words in $Y^{\ast}$ which contains $z_p$ $a$ times and $z_{2p}$ $b$ times. Then for any nonnegative integers $m,n$, we have
\begin{align}
z_p^m \ast z_p^n=\sum\limits_{i=0}^{\min\{m,n\}}\binom{m+n-2i}{m-i}\sum\limits_{w\in G_i^{m+n-2i}}w
\label{Eq:SimPro-MZV}
\end{align}
and
\begin{align}
z_p^m \,\bar{\ast}\, z_p^n=\sum\limits_{i=0}^{\min\{m,n\}}{(-1)}^i\binom{m+n-2i}{m-i}\sum\limits_{w\in G_i^{m+n-2i}}w.
\label{Eq:SimPro-MZSV}
\end{align}
\end{thm}

Note that formula \eqref{Eq:SimPro-MZV} is just \cite[Lemma  1]{Chen}. To study the general stuffle products mentioned above, we need the following recursive formulas, which are similar to \cite[Lemma 2.1]{Lei-Guo-Ma} for shuffle products.

\begin{thm}\label{Thm:Recursive-Formula}
Let $m,n, k_1,\ldots,k_m,l_1,\ldots,l_n$ be positive integers. Then for any integer $j$ with $1\leqslant j\leqslant m$,  we have
\begin{align}
&z_{k_1}z_{k_2}\cdots z_{k_m} \ast z_{l_1}z_{l_2}\cdots z_{l_n}=\sum\limits_{i=0}^n[(z_{k_1}z_{k_2}\cdots z_{k_{j-1}}\ast z_{l_1}z_{l_2}\cdots z_{l_i})z_{k_j}\nonumber\\
&\quad +(z_{k_1}z_{k_2}\cdots z_{k_{j-1}}\ast z_{l_1}z_{l_2}\cdots z_{l_{i-1}})z_{k_j+l_i}](z_{k_{j+1}}\cdots z_{k_{m}}\ast z_{l_{i+1}}\cdots z_{l_{n}}),
\label{Eq:RecurFormula-MZV}
\end{align}
and
\begin{align}
&z_{k_1}z_{k_2}\cdots z_{k_m} \,\bar{\ast}\, z_{l_1}z_{l_2}\cdots z_{l_n}=\sum\limits_{i=0}^n[(z_{k_1}z_{k_2}\cdots z_{k_{j-1}}\,\bar{\ast}\, z_{l_1}z_{l_2}\cdots z_{l_i})z_{k_j}\nonumber\\
&\quad -(z_{k_1}z_{k_2}\cdots z_{k_{j-1}}\,\bar{\ast}\, z_{l_1}z_{l_2}\cdots z_{l_{i-1}})z_{k_j+l_i}](z_{k_{j+1}}\cdots z_{k_{m}}\,\bar{\ast}\, z_{l_{i+1}}\cdots z_{l_{n}}).
\label{Eq:RecurFormula-MZSV}
\end{align}
Here by convention, we let $z_{k_i}\cdots z_{k_{i-1}}=1$ and $z_{k_i}\cdots z_{k_{i-2}}=0$ for any integer $i$.
\end{thm}

Applying the recursive formulas to $z_kz_p^m\ast z_lz_p^n$ and $z_kz_p^m\,\bar{\ast}\, z_lz_p^n$, we get the following stuffle product formulas.

\begin{cor}\label{Cor:Stuffle-1-1}
Let $k,l,p$ be positive integers and let $m,n$ be nonnegative integers. Then we have
\begin{align}
&z_k z_p^m\ast z_l z_p^n=\sum\limits_{i=1}^m z_k (z_p^{i}z_l+z_p^{i-1}z_{p+l})(z_p^{m-i}\ast z_p^n)\nonumber\\
&\quad+\sum\limits_{i=1}^{n}z_l (z_p^{i}z_k+z_p^{i-1}z_{p+k})(z_p^{n-i}\ast z_p^m)+(z_k z_l+z_l z_k+z_{k+l})(z_p^m\ast z_p^n)
\label{Eq:Stu-1-1-MZV}
\end{align}
and
\begin{align}
& z_k z_p^m\,\bar{\ast}\, z_l z_p^n=\sum\limits_{i=1}^m z_k (z_p^{i}z_l-z_p^{i-1}z_{p+l})(z_p^{m-i}\,\bar{\ast}\, z_p^n)\nonumber\\
&\quad+\sum\limits_{i=1}^{n}z_l (z_p^{i}z_k-z_p^{i-1}z_{p+k})(z_p^{n-i}\,\bar{\ast}\, z_p^m)+(z_k z_l+z_l z_k-z_{k+l})(z_p^m\,\bar{\ast}\, z_p^n).
\label{Eq:Stu-1-1-MZSV}
\end{align}
\end{cor}

Note that formula \eqref{Eq:Stu-1-1-MZV} is just \cite[Theorem 3]{Chen}. Applying the recursive formulas to $z_kz_p^m\ast z_{l_1}z_p^{n_1}z_{l_2}z_p^{n_2}$ and $z_kz_p^m\,\bar{\ast}\, z_{l_1}z_p^{n_1}z_{l_2}z_p^{n_2}$, we get the following stuffle product formulas.

\begin{cor}\label{Cor:Stuffle-1-2}
Let $k,l_1,l_2,p$ be positive integers and let $m,n_1,n_2$ be nonnegative integers. Then we have
\begin{align}
&z_k z_p^m \ast z_{l_1} z_p^{n_1} z_{l_2} z_p^{n_2}
=\sum\limits_{0\leqslant i_1\leqslant i_2\leqslant m}z_k(z_p^{i_1} z_{l_1}+z_p^{i_1-1}z_{p+l_1})[(z_p^{i_2-i_1}\ast z_p^{n_1})z_{l_2}\nonumber\\
&\;+(z_p^{i_2-i_1-1}\ast z_p^{n_1})z_{p+l_2}](z_p^{m-i_2}\ast z_p^{n_2})+\sum\limits_{0\leqslant i_1\leqslant {n_1}\atop 0\leqslant i_2\leqslant m}z_{l_1}(z_p^{i_1} z_k+z_p^{i_1-1}z_{p+k})\nonumber\\
&\;\times [(z_p^{i_2}\ast z_p^{n_1-i_1})z_{l_2}+(z_p^{i_2-1}\ast z_p^{n_1-i_1})z_{p+l_2}](z_p^{m-i_2}\ast z_p^{n_2})\nonumber\\
&\;+z_{l_1}z_p^{n_1}z_{k+l_2}(z_{p}^m\ast z_p^{n_2})+\sum\limits_{0\leqslant i\leqslant n_2}z_{l_1}z_p^{n_1}z_{l_2}(z_p^{i} z_k+z_p^{i-1}z_{p+k})(z_p^m\ast z_p^{n_2-i})\nonumber\\
&\;+\sum\limits_{0\leqslant i\leqslant m}z_{k+l_1}[(z_p^{i}\ast z_p^{n_1})z_{l_2}+(z_p^{i-1}\ast z_p^{n_1})z_{p+l_2}](z_p^{m-i}\ast z_p^{n_2})
\label{Eq:Stu-1-2-MZV}
\end{align}
and
\begin{align}
&z_k z_p^m \,\bar{\ast}\, z_{l_1} z_p^{n_1} z_{l_2} z_p^{n_2}
=\sum\limits_{0\leqslant i_1\leqslant i_2\leqslant m}z_k(z_p^{i_1} z_{l_1}-z_p^{i_1-1}z_{p+l_1})[(z_p^{i_2-i_1}\,\bar{\ast}\, z_p^{n_1})z_{l_2}\nonumber\\
&\;-(z_p^{i_2-i_1-1}\,\bar{\ast}\, z_p^{n_1})z_{p+l_2}](z_p^{m-i_2}\,\bar{\ast}\, z_p^{n_2})+\sum\limits_{0\leqslant i_1\leqslant {n_1}\atop 0\leqslant i_2\leqslant m}z_{l_1}(z_p^{i_1} z_k-z_p^{i_1-1}z_{p+k})\nonumber\\
&\;\times [(z_p^{i_2}\,\bar{\ast}\, z_p^{n_1-i_1})z_{l_2}-(z_p^{i_2-1}\,\bar{\ast}\, z_p^{n_1-i_1})z_{p+l_2}](z_p^{m-i_2}\,\bar{\ast}\, z_p^{n_2})\nonumber\\
&\;-z_{l_1}z_p^{n_1}z_{k+l_2}(z_{p}^m\,\bar{\ast}\, z_p^{n_2})+\sum\limits_{0\leqslant i\leqslant n_2}z_{l_1}z_p^{n_1}z_{l_2}(z_p^{i} z_k-z_p^{i-1}z_{p+k})(z_p^m\,\bar{\ast}\, z_p^{n_2-i})\nonumber\\
&\;-\sum\limits_{0\leqslant i\leqslant m}z_{k+l_1}[(z_p^{i}\,\bar{\ast}\, z_p^{n_1})z_{l_2}-(z_p^{i-1}\,\bar{\ast}\, z_p^{n_1})z_{p+l_2}](z_p^{m-i}\,\bar{\ast}\, z_p^{n_2}).
\label{Eq:Stu-1-2-MZSV}
\end{align}
Here we set $z_p^{-1}=0$.
\end{cor}

Similarly, one can use the recursive formula \eqref{Eq:RecurFormula-MZV} to the general product \eqref{Eq:Stu-MZV} and  formula \eqref{Eq:RecurFormula-MZSV} to the general product \eqref{Eq:Stu-MZSV}. We find that all these type products can be expressed by simple products of the forms $z_p^i\ast z_p^j$ for the $\ast$ product, and by simple products of the forms $z_p^i\,\bar{\ast}\,z_p^j$ for the $\bar{\ast}$ product. However the complexity of the formula increases quickly when $m+n$ is larger. Hence we do not write the explicit formulas here.

After getting the formula of the product \eqref{Eq:Stu-MZV}, one can apply the map $Z$ to get a formula of the stuffle product
$$\zeta(k_1,\{p\}^{m_1},\ldots,k_r,\{p\}^{m_r})\zeta(l_1,\{p\}^{n_1},\ldots,l_s,\{p\}^{n_s}),$$
where $k_1,l_1\geqslant 2$ and $\{p\}^m$ stands for a string of $m$ $p$'s. Similarly, from the formula of the product \eqref{Eq:Stu-MZSV}, one can obtain a corresponding formula of the stuffle product of multiple zeta-star values. For example, from Theorem \ref{Thm:Simple-Product}, Corollary \ref{Cor:Stuffle-1-1} and Corollary \ref{Cor:Stuffle-1-2}, we immediately get the following corollary. Note that \eqref{Eq:Stuffle-MZV} is just \cite[Theorem 1]{Chen}.

\begin{cor}
Let $p$ be a positive integer. For nonnegative integers $a,b$, let $A_b^a$ be the set of all strings containing $a$ times $p$ and $b$ times $2p$. Then for any positive integers $k,l,l_1,l_2$ with $k,l,l_1\geqslant 2$ and any nonnegative integers $m,n,n_1,n_2$, we have
\begin{align}
&\zeta(k,\{p\}^m)\zeta(l,\{p\}^n)\nonumber\\
=&\sum\limits_{{1\leqslant i\leqslant m\atop 0\leqslant j\leqslant\min\{m-i,n\}}\atop \mathbf{k}\in A^{m+n-i-2j}_j}\binom{m+n-i-2j}{n-j}\left[\zeta(k,\{p\}^i,l,\mathbf{k})+\zeta(k,\{p\}^{i-1},p+l,\mathbf{k})\right]\nonumber\\
&+\sum\limits_{{1\leqslant i\leqslant n\atop 0\leqslant j\leqslant \min\{m,n-i\}}\atop \mathbf{k}\in A^{m+n-i-2j}_j}\binom{m+n-i-2j}{m-j}\left[\zeta(l,\{p\}^i,k,\mathbf{k})+\zeta(l,\{p\}^{i-1},p+k,\mathbf{k})\right]\nonumber\\
&+\sum\limits_{0\leqslant j\leqslant \min\{m,n\}\atop \mathbf{k}\in A^{m+n-2j}_j}\binom{m+n-2j}{m-j}\left[\zeta(k,l,\mathbf{k})+\zeta(l,k,\mathbf{k})+\zeta(k+l,\mathbf{k})\right],\label{Eq:Stuffle-MZV}\\
&\zeta^{\star}(k,\{p\}^m)\zeta^{\star}(l,\{p\}^n)\nonumber\\
=&\sum\limits_{{1\leqslant i\leqslant m\atop 0\leqslant j\leqslant\min\{m-i,n\}}\atop \mathbf{k}\in A^{m+n-i-2j}_j}(-1)^j\binom{m+n-i-2j}{n-j}\left[\zeta^{\star}(k,\{p\}^i,l,\mathbf{k})-\zeta^{\star}(k,\{p\}^{i-1},p+l,\mathbf{k})\right]\nonumber\\
&+\sum\limits_{{1\leqslant i\leqslant n\atop 0\leqslant j\leqslant \min\{m,n-i\}}\atop \mathbf{k}\in A^{m+n-i-2j}_j}(-1)^j\binom{m+n-i-2j}{m-j}\left[\zeta^{\star}(l,\{p\}^i,k,\mathbf{k})\right.\nonumber\\
&\qquad\qquad\qquad\qquad\qquad\qquad\qquad\qquad\qquad\qquad\left.-\zeta^{\star}(l,\{p\}^{i-1},p+k,\mathbf{k})\right]\nonumber\\
&+\sum\limits_{0\leqslant j\leqslant \min\{m,n\}\atop \mathbf{k}\in A^{m+n-2j}_j}(-1)^j\binom{m+n-2j}{m-j}\left[\zeta^{\star}(k,l,\mathbf{k})+\zeta^{\star}(l,k,\mathbf{k})-\zeta^{\star}(k+l,\mathbf{k})\right]\nonumber
\end{align}
and
\begin{align*}
&\zeta(k,\{p\}^m)\zeta(l_1,\{p\}^{n_1},l_2,\{p\}^{n_2})\\
=&\sum\limits_{{{0\leqslant i_1\leqslant i_2\leqslant m\atop 0\leqslant j_1\leqslant \min\{i_2-i_1,n_1\}}\atop{0\leqslant j_2\leqslant \min\{m-i_2,n_2\}\atop \mathbf{k}\in A^{n_1+i_2-i_1-2j_1}_{j_1}}}\atop \mathbf{l}\in A^{m+n_2-i_2-2j_2}_{j_2}}\binom{n_1+i_2-i_1-2j_1}{n_1-j_1}\binom{m+n_2-i_2-2j_2}{n_2-j_2}\\
&\quad\times \left[\zeta(k,\{p\}^{i_1},l_1,\mathbf{k},l_2,\mathbf{l})+\zeta(k,\{p\}^{i_1-1},p+l_1,\mathbf{k},l_2,\mathbf{l})\right]\\
&+\sum\limits_{{{0\leqslant i_1\leqslant i_2\leqslant m\atop 0\leqslant j_1\leqslant \min\{i_2-i_1-1,n_1\}}\atop{0\leqslant j_2\leqslant \min\{m-i_2,n_2\}\atop \mathbf{k}\in A^{n_1+i_2-i_1-2j_1-1}_{j_1}}}\atop \mathbf{l}\in A^{m+n_2-i_2-2j_2}_{j_2}}\binom{n_1+i_2-i_1-2j_1-1}{n_1-j_1}\binom{m+n_2-i_2-2j_2}{n_2-j_2}\\
&\quad\times \left[\zeta(k,\{p\}^{i_1},l_1,\mathbf{k},p+l_2,\mathbf{l})+\zeta(k,\{p\}^{i_1-1},p+l_1,\mathbf{k},p+l_2,\mathbf{l})\right]\\
&+\sum\limits_{{{0\leqslant i_1\leqslant n_1,0\leqslant i_2\leqslant m\atop 0\leqslant j_1\leqslant \min\{i_2,n_1-i_1\}}\atop{0\leqslant j_2\leqslant \min\{m-i_2,n_2\}\atop \mathbf{k}\in A^{n_1+i_2-i_1-2j_1}_{j_1}}}\atop \mathbf{l}\in A^{m+n_2-i_2-2j_2}_{j_2}}\binom{n_1+i_2-i_1-2j_1}{i_2-j_1}\binom{m+n_2-i_2-2j_2}{n_2-j_2}\\
&\quad\times\left[\zeta(l_1,\{p\}^{i_1},k,\mathbf{k},l_2,\mathbf{l})+\zeta(l_1,\{p\}^{i_1-1},p+k,\mathbf{k},l_2,\mathbf{l})\right]\\
&+\sum\limits_{{{0\leqslant i_1\leqslant n_1,0\leqslant i_2\leqslant m\atop 0\leqslant j_1\leqslant \min\{i_2-1,n_1-i_1\}}\atop{0\leqslant j_2\leqslant \min\{m-i_2,n_2\}\atop \mathbf{k}\in A^{n_1+i_2-i_1-2j_1-1}_{j_1}}}\atop \mathbf{l}\in A^{m+n_2-i_2-2j_2}_{j_2}}\binom{n_1+i_2-i_1-2j_1-1}{i_2-j_1-1}\binom{m+n_2-i_2-2j_2}{n_2-j_2}\\
&\quad\times\left[\zeta(l_1,\{p\}^{i_1},k,\mathbf{k},p+l_2,\mathbf{l})+\zeta(l_1,\{p\}^{i_1-1},p+k,\mathbf{k},p+l_2,\mathbf{l})\right]\\
&+\sum\limits_{0\leqslant j\leqslant \min\{m,n_2\}\atop \mathbf{k}\in A^{m+n_2-2j}_j}\binom{m+n_2-2j}{m-j}\zeta(l_1,\{p\}^{n_1},k+l_2,\mathbf{k})\\
&+\sum\limits_{{0\leqslant i\leqslant n_2\atop 0\leqslant j\leqslant \min\{m,n_2-i\}}\atop \mathbf{k}\in A^{m+n_2-i-2j}_j}\binom{m+n_2-i-2j}{m-j}\\
&\quad\times\left[\zeta(l_1,\{p\}^{n_1},l_2,\{p\}^{i},k,\mathbf{k})+\zeta(l_1,\{p\}^{n_1},l_2,\{p\}^{i-1},p+k,\mathbf{k})\right]\\
&+\sum\limits_{{{0\leqslant i\leqslant m\atop 0\leqslant j_1\leqslant \min\{i,n_1\}}\atop{0\leqslant j_2\leqslant \min\{m-i,n_2\}\atop\mathbf{k}\in A^{n_1+i-2j_1}_{j_1}}}\atop \mathbf{l}\in A^{m+n_2-i-2j_2}_{j_2}}\binom{n_1+i-2j_1}{n_1-j_1}\binom{m+n_2-i-2j_2}{n_2-j_2}\zeta(k+l_1,\mathbf{k},l_2,\mathbf{l})\\
&+\sum\limits_{{{0\leqslant i\leqslant m\atop 0\leqslant j_1\leqslant \min\{i-1,n_1\}}\atop{0\leqslant j_2\leqslant \min\{m-i,n_2\}\atop\mathbf{k}\in A^{n_1+i-2j_1-1}_{j_1}}}\atop \mathbf{l}\in A^{m+n_2-i-2j_2}_{j_2}}\binom{n_1+i-2j_1-1}{n_1-j_1}\binom{m+n_2-i-2j_2}{n_2-j_2}\zeta(k+l_1,\mathbf{k},p+l_2,\mathbf{l}),\\
&\zeta^{\star}(k,\{p\}^m)\zeta^{\star}(l_1,\{p\}^{n_1},l_2,\{p\}^{n_2})\\
=&\sum\limits_{{{0\leqslant i_1\leqslant i_2\leqslant m\atop 0\leqslant j_1\leqslant \min\{i_2-i_1,n_1\}}\atop{0\leqslant j_2\leqslant \min\{m-i_2,n_2\}\atop \mathbf{k}\in A^{n_1+i_2-i_1-2j_1}_{j_1}}}\atop \mathbf{l}\in A^{m+n_2-i_2-2j_2}_{j_2}}(-1)^{j_1+j_2}\binom{n_1+i_2-i_1-2j_1}{n_1-j_1}\binom{m+n_2-i_2-2j_2}{n_2-j_2}\\
&\quad\times \left[\zeta^{\star}(k,\{p\}^{i_1},l_1,\mathbf{k},l_2,\mathbf{l})-\zeta^{\star}(k,\{p\}^{i_1-1},p+l_1,\mathbf{k},l_2,\mathbf{l})\right]\\
&-\sum\limits_{{{0\leqslant i_1\leqslant i_2\leqslant m\atop 0\leqslant j_1\leqslant \min\{i_2-i_1-1,n_1\}}\atop{0\leqslant j_2\leqslant \min\{m-i_2,n_2\}\atop \mathbf{k}\in A^{n_1+i_2-i_1-2j_1-1}_{j_1}}}\atop \mathbf{l}\in A^{m+n_2-i_2-2j_2}_{j_2}}(-1)^{j_1+j_2}\binom{n_1+i_2-i_1-2j_1-1}{n_1-j_1}\binom{m+n_2-i_2-2j_2}{n_2-j_2}\\
&\quad\times \left[\zeta^{\star}(k,\{p\}^{i_1},l_1,\mathbf{k},p+l_2,\mathbf{l})-\zeta^{\star}(k,\{p\}^{i_1-1},p+l_1,\mathbf{k},p+l_2,\mathbf{l})\right]\\
&+\sum\limits_{{{0\leqslant i_1\leqslant n_1,0\leqslant i_2\leqslant m\atop 0\leqslant j_1\leqslant \min\{i_2,n_1-i_1\}}\atop{0\leqslant j_2\leqslant \min\{m-i_2,n_2\}\atop \mathbf{k}\in A^{n_1+i_2-i_1-2j_1}_{j_1}}}\atop \mathbf{l}\in A^{m+n_2-i_2-2j_2}_{j_2}}(-1)^{j_1+j_2}\binom{n_1+i_2-i_1-2j_1}{i_2-j_1}\binom{m+n_2-i_2-2j_2}{n_2-j_2}\\
&\quad\times\left[\zeta^{\star}(l_1,\{p\}^{i_1},k,\mathbf{k},l_2,\mathbf{l})-\zeta^{\star}(l_1,\{p\}^{i_1-1},p+k,\mathbf{k},l_2,\mathbf{l})\right]\\
&-\sum\limits_{{{0\leqslant i_1\leqslant n_1,0\leqslant i_2\leqslant m\atop 0\leqslant j_1\leqslant \min\{i_2-1,n_1-i_1\}}\atop{0\leqslant j_2\leqslant \min\{m-i_2,n_2\}\atop \mathbf{k}\in A^{n_1+i_2-i_1-2j_1-1}_{j_1}}}\atop \mathbf{l}\in A^{m+n_2-i_2-2j_2}_{j_2}}(-1)^{j_1+j_2}\binom{n_1+i_2-i_1-2j_1-1}{i_2-j_1-1}\binom{m+n_2-i_2-2j_2}{n_2-j_2}\\
&\quad\times\left[\zeta^{\star}(l_1,\{p\}^{i_1},k,\mathbf{k},p+l_2,\mathbf{l})-\zeta^{\star}(l_1,\{p\}^{i_1-1},p+k,\mathbf{k},p+l_2,\mathbf{l})\right]\\
&-\sum\limits_{0\leqslant j\leqslant \min\{m,n_2\}\atop \mathbf{k}\in A^{m+n_2-2j}_j}(-1)^{j}\binom{m+n_2-2j}{m-j}\zeta^{\star}(l_1,\{p\}^{n_1},k+l_2,\mathbf{k})\\
&+\sum\limits_{{0\leqslant i\leqslant n_2\atop 0\leqslant j\leqslant \min\{m,n_2-i\}}\atop \mathbf{k}\in A^{m+n_2-i-2j}_j}(-1)^{j}\binom{m+n_2-i-2j}{m-j}\\
&\quad\times\left[\zeta^{\star}(l_1,\{p\}^{n_1},l_2,\{p\}^{i},k,\mathbf{k})-\zeta^{\star}(l_1,\{p\}^{n_1},l_2,\{p\}^{i-1},p+k,\mathbf{k})\right]\\
&-\sum\limits_{{{0\leqslant i\leqslant m\atop 0\leqslant j_1\leqslant \min\{i,n_1\}}\atop{0\leqslant j_2\leqslant \min\{m-i,n_2\}\atop\mathbf{k}\in A^{n_1+i-2j_1}_{j_1}}}\atop \mathbf{l}\in A^{m+n_2-i-2j_2}_{j_2}}(-1)^{j_1+j_2}\binom{n_1+i-2j_1}{n_1-j_1}\binom{m+n_2-i-2j_2}{n_2-j_2}\zeta^{\star}(k+l_1,\mathbf{k},l_2,\mathbf{l})\\
&+\sum\limits_{{{0\leqslant i\leqslant m\atop 0\leqslant j_1\leqslant \min\{i-1,n_1\}}\atop{0\leqslant j_2\leqslant \min\{m-i,n_2\}\atop\mathbf{k}\in A^{n_1+i-2j_1-1}_{j_1}}}\atop \mathbf{l}\in A^{m+n_2-i-2j_2}_{j_2}}(-1)^{j_1+j_2}\binom{n_1+i-2j_1-1}{n_1-j_1}\binom{m+n_2-i-2j_2}{n_2-j_2}\\
&\qquad\qquad\qquad\qquad\qquad\qquad\qquad\qquad\times\zeta^{\star}(k+l_1,\mathbf{k},p+l_2,\mathbf{l}).
\end{align*}
\end{cor}

 Also, if we combine the stuffle product formulas obtained here together with shuffle product formulas obtained in  \cite{Eie-Wei} or \cite{Li-Qin}, we can get some double shuffle relations of multiple zeta values for some special types. Here we leave the explicit formulas to the reader.

Section \ref{Sec:Proof} contains proofs of Theorem \ref{Thm:Simple-Product}, Theorem \ref{Thm:Recursive-Formula}, Corollary \ref{Cor:Stuffle-1-1} and Corollary \ref{Cor:Stuffle-1-2}. In Section \ref{Sec:GenProd}, we describe how to apply our recursive formulas to express the products \eqref{Eq:Stu-MZV} and \eqref{Eq:Stu-MZSV} by simple products, and take the cases $r=s=2$ as an example.


\section{Proofs} \label{Sec:Proof}

In this section, we give proofs of the results mentioned in Section \ref{Sec:Intro}.

\subsection{Proof of Theorem \ref{Thm:Simple-Product}}

For the simple product $z_p^m\ast z_p^n$, Chen used induction on $m+n$ to prove formula \eqref{Eq:SimPro-MZV} in \cite{Chen}. Here we give another proof, which seems simple and intuitive. We use the combinatorial description of stuffle product $\ast$. By the definition of stuffle product $\ast$, we have
$$z_{k_1}\cdots z_{k_n}\ast z_{k_{n+1}}\cdots z_{k_{n+m}}=\sum\limits_{i=0}^{\min\{n,m\}}\sum\limits_{\sigma\in\mathfrak{S}_{n,m,i}}z_{\sigma^{-1}(1)}z_{\sigma^{-1}(2)}\cdots z_{\sigma^{-1}(n+m-i)},$$
where $\mathfrak{S}_{n,m,i}$ is the set of all surjections $\sigma: \{1,2,\ldots,n+m\}\twoheadrightarrow \{1,2,\ldots,n+m-i\}$ with conditions
$$\sigma(1)<\sigma(2)<\cdots<\sigma(n),\quad
\sigma(n+1)<\sigma(n+2)<\cdots<\sigma(n+m),$$
and for any $\sigma\in \mathfrak{S}_{n,m,i}$ and any $j\in \{1,2,\ldots,n+m-i\}$,
$$z_{\sigma^{-1}(j)}=\begin{cases}
z_{k}, & \text{if\;} \sigma^{-1}(j)=\{k\},\\
z_{k+l}, & \text{if\;} \sigma^{-1}(j)=\{k,l\}.
\end{cases}$$

Now if we write
$$z_p^m\ast z_p^n=\sum\limits_{w\in Y^{\ast}}c_ww,$$
then we only need to determine the coefficients $c_w$. By the combinatorial description of the product $\ast$ described above, it is easy to see that if $c_w\neq 0$, then all the letters in the word $w$ must be $z_p$ or $z_{2p}$, and the number of $z_{2p}$'s is less than or equal to $\min\{m,n\}$. Hence we can write
$$z_p^m \ast z_p^n=\sum\limits_{i=0}^{\min\{m,n\}}\sum\limits_{w\in G_i^{m+n-2i}}c_i^w w.$$

For any $w\in G_i^{m+n-2i}$, we look for the possibility to get $w$. There are $m+n-2i$ times $z_p$ in $w$, in which $m-i$ times $z_p$ are from $z_p^m$ and others are from $z_p^n$. Hence the possibility is $\binom{m+n-2i}{m-i}$. In other words, for any $w\in G_i^{m+n-2i}$, we have
$$c_i^w=\binom{m+n-2i}{m-i}.$$
Then we get \eqref{Eq:SimPro-MZV}.

Similarly, to prove formula \eqref{Eq:SimPro-MZSV} for the simple product $z_p^m\,\bar{\ast}\,z_p^n$, we use the combinatorial description of stuffle product $\bar{\ast}$. we have
$$z_{k_1}\cdots z_{k_n}\,\bar{\ast}\, z_{k_n+1}\cdots z_{k_{n+m}}=\sum\limits_{i=0}^{\min\{n,m\}}\sum\limits_{\sigma\in\mathfrak{S}_{n,m,i}}z_{\sigma^{-1}(1)}^{\star}z_{\sigma^{-1}(2)}^{\star}\cdots z_{\sigma^{-1}(n+m-i)}^{\star},$$
where for any $\sigma\in \mathfrak{S}_{n,m,i}$ and any $j\in \{1,2,\ldots,n+m-i\}$,
$$z_{\sigma^{-1}(j)}^{\star}=\begin{cases}
z_{k}, & \text{if\;} \sigma^{-1}(j)=\{k\},\\
-z_{k+l}, & \text{if\;} \sigma^{-1}(j)=\{k,l\}.
\end{cases}$$

By the combinatorial description, we can write
$$z_p^m\, \bar{\ast}\, z_p^n=\sum\limits_{i=0}^{\min\{m,n\}}\sum\limits_{w\in G_i^{m+n-2i}}d_i^w w.$$
In this case, for any $w\in G_i^{m+n-2i}$, the possibility to get $w$ is also $\binom{m+n-2i}{m-i}$. While every letter $z_{2p}$ in $w$ will bring a negative sign. Hence  for any $w\in G_i^{m+n-2i}$, we have
$$d_i^w=(-1)^i\binom{m+n-2i}{m-i}.$$
Then we get \eqref{Eq:SimPro-MZSV}.\qed

For more details, see \cite{Li-Qin}, which uses the same idea to deal with shuffle products.

\subsection{Proof of Theorem \ref{Thm:Recursive-Formula}}

We only prove the recursive formula \eqref{Eq:RecurFormula-MZV}. The recursive formula \eqref{Eq:RecurFormula-MZSV} can be proven in a similar way. To get the result, one can use induction on $m+n$, while here we use the combinatorial description of the stuffle product.

The left-hand side of \eqref{Eq:RecurFormula-MZV} is a sum over surjections
$$\sigma:\{1,2,\ldots,m+n\}\twoheadrightarrow \{1,2,\ldots,m+n-p\},$$
which are strictly increasing on $\{1,\ldots,m\}$ and on $\{m+1,\ldots,m+n\}$ separately, and $p$ is an integer with the condition $0\leqslant p\leqslant \min\{m,n\}$.
Now  for  a fixed $j\in \{1,2,\dots,m\}$, there are two different types of such $\sigma$'s, one satisfies the condition $\sigma(i)< \sigma(j)< \sigma(i+1)$ for some $i\in\{m,m+1,\ldots,m+n\}$, and the other one satisfies the condition $\sigma(j)=\sigma(i)$ for some $i\in \{m+1,m+2,\dots,m+n\}$. Summing over the first type of $\sigma$'s we get the term
$$\sum\limits_{i=0}^n(z_{k_1}z_{k_2}\cdots z_{k_{j-1}}\ast z_{l_1}z_{l_2}\cdots z_{l_i})z_{k_j}(z_{k_{j+1}}\cdots z_{k_{m}}\ast z_{l_{i+1}}\cdots z_{l_{n}}).$$
And summing over the second type of $\sigma$'s we obtain the term
$$\sum\limits_{i=1}^n(z_{k_1}z_{k_2}\cdots z_{k_{j-1}}\ast z_{l_1}z_{l_2}\cdots z_{l_{i-1}})z_{k_j+l_i}(z_{k_{j+1}}\cdots z_{k_{m}}\ast z_{l_{i+1}}\cdots z_{l_{n}}).$$
Then we get the recursive formula \eqref{Eq:RecurFormula-MZV}.\qed

\subsection{Proof of Corollary \ref{Cor:Stuffle-1-1}}

Here we prove formula \eqref{Eq:Stu-1-1-MZV}, and one can prove formula \eqref{Eq:Stu-1-1-MZSV} in a similar way. We remark that the  strategy of our proof is slightly different from that of \cite[Theorem 3]{Chen}. We apply the recursive formula \eqref{Eq:RecurFormula-MZV} with $j=1$ to get
\begin{align*}
z_kz_p^m\ast z_lz_p^n=&z_kz_l(z_p^m\ast z_p^n)+z_{k+l}(z_p^m\ast z_p^n)\\
&\qquad+\sum\limits_{i=0}^n(z_lz_p^iz_k+z_lz_p^{i-1}z_{p+k})(z_p^m\ast z_p^{n-i}),
\end{align*}
which is just \eqref{Eq:Stu-1-1-MZV}.\qed

\subsection{Prove of Corollary \ref{Cor:Stuffle-1-2}}

Here we prove formula \eqref{Eq:Stu-1-2-MZV}, and one can prove formula \eqref{Eq:Stu-1-2-MZSV} in a similar way. Applying the recursive formula \eqref{Eq:RecurFormula-MZV} with $j=1$, we get
\begin{align}
&z_k z_p^m \ast z_{l_1} z_p^{n_1} z_{l_2} z_p^{n_2}=z_k(z_p^m\ast z_{l_1}z_p^{n_1}z_{l_2}z_p^{n_2})+z_{k+l_1}(z_p^m\ast z_p^{n_1} z_{l_2} z_p^{n_2})\nonumber\\
&\qquad +\sum\limits_{i=0}^{n_1}(z_{l_1}z_p^iz_{k}+z_{l_1}z_p^{i-1}z_{p+k})(z_p^m\ast z_p^{n_1-i}z_{l_2}z_p^{n_2})+z_{l_1}z_p^{n_1}z_{k+l_2}(z_p^m\ast z_p^{n_2})\nonumber\\
&\qquad +\sum\limits_{i=0}^{n_2}(z_{l_1}z_p^{n_1}z_{l_2}z_p^iz_k+z_{l_1}z_p^{n_1}z_{l_2}z_p^{i-1}z_{p+k})(z_p^m\ast z_p^{n_2-i}).
\label{Eq:Stu-1-2-Step1}
\end{align}
Then we have to compute the products $z_p^m\ast z_{l_1}z_p^{n_1}z_{l_2}z_p^{n_2}$ and $z_p^m\ast z_p^{n_1-i}z_{l_2}z_p^{n_2}$. For the second one, we have the following general result.

\begin{lem}
For positive integers $l,p$ and nonnegative integers $m,n_1,n_2$, we have
\begin{align}
z_p^m\ast z_p^{n_1} z_l z_p^{n_2}=\sum\limits_{i=0}^m[(z_p^{i}\ast z_p^{n_1})z_l+(z_p^{i-1}\ast z_p^{n_1})z_{p+l}](z_p^{m-i}\ast z_p^{n_2}).
\label{Eq:Stu-0-2-MZV}
\end{align}
\end{lem}

\proof We get \eqref{Eq:Stu-0-2-MZV} by applying the recursive formula \eqref{Eq:RecurFormula-MZV} with $j=n_1+1$. \qed

And for the first one, we have the following result.

\begin{lem}
For positive integers $p,l_1,l_2$ and nonnegative integers $m,n_1,n_2$, we have
\begin{align}
&z_p^m\ast z_{l_1} z_p^{n_1} z_{l_2} z_p^{n_2}=\sum\limits_{0\leqslant i_1\leqslant i_2\leqslant m}(z_p^{i_1} z_{l_1}+z_p^{i_1-1}z_{p+{l_1}})\nonumber\\
&\qquad \times [(z_p^{i_2-i_1}\ast z_p^{n_1})z_{l_2}+(z_p^{ i_2- i_1-1}\ast z_p^{n_1})z_{p+{l_2}}](z_p^{m-i_2}\ast z_p^{n_2}).
\label{Eq:Stu-0-22-MZV}
\end{align}
\end{lem}

\proof Applying the recursive formula \eqref{Eq:RecurFormula-MZV} with $j=1$, we get
\begin{align*}
z_p^m\ast z_{l_1} z_p^{n_1} z_{l_2} z_p^{n_2}=\sum\limits_{i=0}^m(z_p^{i} z_{l_1}+z_p^{i-1}z_{p+{l_1}})(z_p^{m-i}\ast z_p^{n_1} z_{l_2} z_p^{n_2}).	
\end{align*}
Then using \eqref{Eq:Stu-0-2-MZV}, we have
\begin{align*}
&z_p^m\ast z_{l_1} z_p^{n_1} z_{l_2} z_p^{n_2}=\sum\limits_{i=0}^m\sum\limits_{t=0}^{m-i}(z_p^{i} z_{l_1}+z_p^{i-1}z_{p+{l_1}})\\
&\qquad \times [(z_p^t\ast z_p^{n_1})z_{l_2}+(z_p^{t-1}\ast z_p^{n_1})z_{p+l_2}](z_p^{m-i-t}\ast z_p^{n_2}).
\end{align*}
Let $i_1=i$ and $i_2=i+t$ in the summation of the right-hand side of the above equation, we get the result. \qed

Now formula \eqref{Eq:Stu-1-2-MZV} follows from \eqref{Eq:Stu-1-2-Step1}, \eqref{Eq:Stu-0-2-MZV} and \eqref{Eq:Stu-0-22-MZV}. And we finish the proof of Corollary \ref{Cor:Stuffle-1-2}.\qed

\subsection{Remark}

In fact, one  can get Corollary \ref{Cor:Stuffle-1-1} and Corollary \ref{Cor:Stuffle-1-2} in the same way as the proofs of Theorem \ref{Thm:Simple-Product} and Theorem \ref{Thm:Recursive-Formula} given above. For example, $z_k z_p^m\ast z_l z_p^n$ is a sum over surjections
$$\sigma:\{1,2,\ldots,m+n+2\}\twoheadrightarrow\{1,2,\ldots,m+n+2-j\},$$
which  are strictly increasing on $\{1,\ldots,m+1\}$ and on $\{m+2,\ldots,m+n+2\}$ respectively. We can separate  such $\sigma$'s into three different types according to
\begin{align*}
&\sigma(1)<\sigma(m+2),\\
&\sigma(m+2)<\sigma(1),\\
\text{or}\quad &\sigma(1)=\sigma(m+2),
\end{align*}
which correspond to the three sums in the right-hand side of \eqref{Eq:Stu-1-1-MZV}. Then we get formula \eqref{Eq:Stu-1-1-MZV}.

But as $r+s$ increases in \eqref{Eq:Stu-MZV},  the possibilities of different types of $\sigma$'s similar to above will increase rapidly. Hence in the following section we discuss the general products by using the recursive formulas.

\section{The general products}\label{Sec:GenProd}

In this section, we discuss the general products \eqref{Eq:Stu-MZV} and \eqref{Eq:Stu-MZSV}. Here we only treat the $\ast$ product. Let $r,s,k_1,\ldots,k_r,l_1,\ldots,l_s,p$ be positive integers and let $m_1,\ldots,m_r$, $n_1,\ldots,n_s$ be nonnegative integers. Applying the recursive formula \eqref{Eq:RecurFormula-MZV} with $j=1$ to the product \eqref{Eq:Stu-MZV}, we find \eqref{Eq:Stu-MZV} can be expressed by products
\begin{itemize}
  \item[($T_1$)] $z_p^{m_1}z_{k_2}z_p^{m_2}\cdots z_{k_r}z_p^{m_r}\ast z_{l_1}z_{p}^{n_1}z_{l_2}z_p^{n_2}\cdots z_{l_s}z_{p}^{n_s}$,
  \item[($T_2$)] $z_p^{m_1}z_{k_2}z_p^{m_2}\cdots z_{k_r}z_p^{m_r}\ast z_{p}^{i_1}z_{l_2}z_p^{n_2}\cdots z_{l_s}z_{p}^{n_s}$,\quad ($0\leqslant i_1\leqslant n_1$),
  \item[($T_3$)] $z_p^{m_1}z_{k_2}z_p^{m_2}\cdots z_{k_r}z_p^{m_r}\ast z_p^{i_2}z_{l_3}\cdots z_{l_s}z_{p}^{n_s}$,\quad ($0\leqslant i_2\leqslant n_2$),
  \item[] $\qquad\qquad\qquad  \vdots$
  \item[($T_{s}$)] $z_p^{m_1}z_{k_2}z_p^{m_2}\cdots z_{k_r}z_p^{m_r}\ast z_{p}^{i_{s-1}}z_{l_{s}}z_p^{n_s}$,\quad ($0\leqslant i_{s-1}\leqslant n_{s-1}$),
  \item[($T_{s+1}$)] $z_p^{m_1}z_{k_2}z_p^{m_2}\cdots z_{k_r}z_p^{m_r}\ast z_{p}^{i_s}$,\quad ($0\leqslant i_s\leqslant n_s$).
\end{itemize}
We may assume that $r>1$. For $i=1,2,\ldots,s$, using the recursive formula \eqref{Eq:RecurFormula-MZV} with $j=m_1+1$ to the product ($T_i$), we find that ($T_i$) can be expressed by products of forms ($T_{i+1}$), $\ldots$, ($T_{s+1}$). For the product ($T_{s+1}$), we also apply the recursive formula \eqref{Eq:RecurFormula-MZV} with $j=m_1+1$ to get
\begin{align*}
&z_p^{m_1}z_{k_2}z_p^{m_2}\cdots z_{k_r}z_p^{m_r}\ast z_{p}^{i_s}\\
=&\sum\limits_{i=0}^{i_s}[(z_p^{m_1}\ast z_p^{i})z_{k_2}+(z_p^{m_1}\ast z_p^{i-1})z_{p+k_2}](z_p^{m_2}z_{k_3}\cdots z_{k_r}z_p^{m_r}\ast z_p^{i_s-i}).
\end{align*}
Hence $(T_{s+1})$ can be expressed by products of forms $(T_{s+1})$ with smaller $r$.
Then we see that we can express the product \eqref{Eq:Stu-MZV} by simple products of the form $z_p^i\ast z_p^j$.

Let's take the case $r=s=2$, that is $z_{k_1}z_p^{m_1}z_{k_2}z_p^{m_2}\ast z_{l_1}z_p^{n_1}z_{l_2}z_p^{n_2}$, as an example. Applying the recursive formula \eqref{Eq:RecurFormula-MZV} with $j=1$, we get
\begin{align}
&z_{k_1}z_p^{m_1}z_{k_2}z_p^{m_2}\ast z_{l_1}z_p^{n_1}z_{l_2}z_p^{n_2}\nonumber\\
=&z_{k_1}(z_p^{m_1}z_{k_2}z_p^{m_2}\ast z_{l_1}z_p^{n_1}z_{l_2}z_p^{n_2})
+z_{k_1+l_1}(z_p^{m_1}z_{k_2}z_p^{m_2}\ast z_p^{n_1}z_{l_2}z_p^{n_2})\nonumber\\
&+\sum\limits_{i=0}^{n_1}(z_{l_1}z_p^iz_{k_1}+z_{l_1}z_p^{i-1}z_{p+k_1})(z_p^{m_1}z_{k_2}z_p^{m_2}\ast z_p^{n_1-i}z_{l_2}z_p^{n_2})\nonumber\\
&+z_{l_1}z_p^{n_1}z_{k_1+l_2}(z_p^{m_1}z_{k_2}z_p^{m_2}\ast z_p^{n_2})\nonumber\\
&+\sum\limits_{i=0}^{n_2}(z_{l_1}z_{p}^{n_1}z_{l_2}z_p^iz_{k_1}+z_{l_1}z_{p}^{n_1}z_{l_2}z_p^{i-1}z_{p+k_1})(z_p^{m_1}z_{k_2}z_p^{m_2}\ast z_p^{n_2-i}).
\label{Eq:Stu-2-2-MZV-Step1}
\end{align}
Hence we need to compute products of the forms
\begin{itemize}
  \item[(i)] $z_p^{m_1}z_{k}z_p^{m_2}\ast z_{l_1}z_p^{n_1}z_{l_2}z_p^{n_2}$,
  \item[(ii)] $z_p^{m_1}z_{k}z_p^{m_2}\ast z_p^{n_1}z_{l}z_p^{n_2}$,
  \item[(iii)] $z_p^{m_1}z_{k}z_p^{m_2}\ast z_p^{n}$.
\end{itemize}

Applying the recursive formula \eqref{Eq:RecurFormula-MZV} with $j=1$ to the product (i), we get
\begin{align*}
&z_p^{m_1}z_{k}z_p^{m_2}\ast z_{l_1}z_p^{n_1}z_{l_2}z_p^{n_2}\\
=&\sum\limits_{i=0}^{m_1}(z_p^iz_{l_1}+z_p^{i-1}z_{p+l_1})(z_p^{m_1-i}z_kz_p^{m_2}\ast z_p^{n_1}z_{l_2}z_p^{n_2})+z_p^{m_1}z_{k+l_1}(z_p^{m_2}\ast z_p^{n_1}z_{l_2}z_p^{n_2})\\
&+\sum\limits_{i=0}^{m_2}(z_p^{m_1}z_kz_p^iz_{l_1}+z_p^{m_1}z_kz_p^{i-1}z_{p+l_1})(z_p^{m_2-i}\ast z_p^{n_1}z_{l_2}z_p^{n_2}).
\end{align*}
Hence the product (i) can be expressed by products of the forms (ii) and (iii). For the product (ii), applying the recursive formula \eqref{Eq:RecurFormula-MZV} with $j=m_1+1$, we get
\begin{align*}
&z_p^{m_1}z_{k}z_p^{m_2}\ast z_p^{n_1}z_{l}z_p^{n_2}\\
=&\sum\limits_{i=0}^{n_1}[(z_p^{m_1}\ast z_p^i)z_k+(z_p^{m_1}\ast z_p^{i-1})z_{p+k}](z_p^{m_2}\ast z_p^{n_1-i}z_lz_p^{n_2})\\
&+(z_p^{m_1}\ast z_p^{n_1})z_{k+l}(z_p^{m_2}\ast z_p^{n_2})\\
&+\sum\limits_{i=0}^{n_2}[(z_p^{m_1}\ast z_p^{n_1}z_lz_p^i)z_k+(z_p^{m_1}\ast z_p^{n_1}z_{l}z_p^{i-1})z_{p+k}](z_p^{m_2}\ast z_p^{n_2-i}).
\end{align*}
Then the product (ii) can be expressed by products of the form (iii). Finally, applying the recursive formula \eqref{Eq:RecurFormula-MZV} with $j=m_1+1$ to the product (iii), we get
\begin{align*}
z_p^{m_1}z_{k}z_p^{m_2}\ast z_p^{n}=\sum\limits_{i=0}^n[(z_p^{m_1}\ast z_p^i)z_k+(z_p^{m_1}\ast z_p^{i-1})z_{p+k}](z_p^{m_2}\ast z_{p}^{n-i}),
\end{align*}
which is expressed by simple products.

Putting all these results into \eqref{Eq:Stu-2-2-MZV-Step1}, we can finally obtain a formula which writes the product $z_{k_1}z_p^{m_1}z_{k_2}z_p^{m_2}\ast z_{l_1}z_p^{n_1}z_{l_2}z_p^{n_2}$ by simple products of the form $z_p^i\ast z_p^j$. While there are too many terms to write, so we omit this formula.


\begin{thebibliography}{99}

\bibitem{Chen} K.-W. Chen, Applications of stuffle product of multiple zeta values, \textit{J. Number Theory} \textbf{153} (2015), 107-116.

\bibitem{Eie-Wei} M. Eie and C.-S. Wei, Generalization of Euler decomposition and their applications, \textit{J. Number Theory} \textbf{133} (2013), 2475-2495.

\bibitem{Hoffman} M. E. Hoffman, The algebra of multiple harmonic series, \textit{J. Algebra} \textbf{194} (2) (1997), 477-495.

\bibitem{Lei-Guo-Ma} P. Lei, L. Guo and B. Ma, Applications of shuffle product to restricted decomposition formulas for multiple zeta values, \textit{J. Number Theory} \textbf{144} (2014), 219-233.

\bibitem{Li-Qin} Z. Li and C. Qin, Shuffle product formulas of multiple zeta values, \textit{J. Number Theory} \textbf{171} (2017), 79-111.

\bibitem{Shuichi Muneta} S. Muneta, Algebraic setup of non-strict multiple zeta values, \textit{Acta Arith.} \textbf{136} (2009), 7-18.

\bibitem{Zagier} D. Zagier, Values of zeta functions and their applications, in First European Congress of Mathematics, Vol. II (Paris, 1992), \textit{Progress in Math.}, Vol. \textbf{120} (Birkh\"{a}user, 1994), pp. 497-512.


\end{thebibliography}
\end{document}